\documentclass[11pt]{article}

\input{epsf}

\usepackage{amsfonts}
\usepackage{amsmath}
\usepackage{amssymb}

\newcommand{\qed}{\hbox{\rule{6pt}{6pt}}}

\newcommand{\Z}{\mathbb{Z}}

\newcommand{\rt}{\vartriangleright}

\newtheorem{theorem}{Theorem}[section]
\newtheorem{corollary}[theorem]{Corollary}

\date{}

\begin{document}

\title{ Automorphism groups of Quandles  }

\author{
Mohamed Elhamdadi 
\\ University of South Florida\footnote{Email: \texttt{emohamed@math.usf.edu}}
\and 
Jennifer MacQuarrie
\\ University of South Florida\footnote{Email: \texttt{jmacquar@mail.usf.edu}}
\and
Ricardo Restrepo\\
Georgia Institute of Technology\footnote{Email: \texttt{ricrest@gmail.com}}
}

\maketitle

\vspace{10mm}

\begin{abstract}

We prove that the automorphism group of the dihedral quandle with $n$ elements is isomorphic to the affine group of the integers mod $n$, and also obtain the inner automorphism group of this quandle. 
 In \cite{HN}, automorphism groups of quandles (up to isomorphisms) of order less than or equal to 5 were given.  
 With the help of the software Maple, we compute the inner and automorphism groups of all $seventy$ $three$ quandles of order six listed in the appendix of  \cite{CKS}.   Since computations of automorphisms of quandles relates to the problem of classification of quandles, we also describe an algorithm implemented in C  for computing all quandles (up to isomorphism) of order less than or equal to nine.  
 
\end{abstract}



\section{Introduction}

Quandles  and racks are algebraic structures whose axiomatization comes from Reidemeister moves in knot theory.  The earliest known work on racks is contained within 1959 correspondence between John Conway and Gavin Wraith who studied racks in the context of  the conjugation operation in a group.  Around 1982,  the notion of a quandle was introduced independently by Joyce \cite{Joyce} and Matveev \cite{Matveev}.  They used it to construct representations of the braid groups.  Joyce and Matveev associated to each knot a quandle that determines the knot up to isotopy and mirror image.  Since then quandles and racks have been investigated by topologists  in order 
to construct knot and link invariants and their higher analogues (see for 
example \cite{CKS} and references therein).

In this paper, we prove that the automorphism group of the dihedral quandle with $n$ elements is isomorphic to the affine group of the integers mod $n$.
  In  \cite{HN}, Ho and Nelson gave the list of quandles (up to isomorphism) of orders $n=3, n=4$ and $n=5$ and determined their automorphism groups. In this paper, with the help of the software Maple, we extend their results by computing the inner and automorphism groups of all $seventy$ $three$ quandles of order six listed in the appendix of  \cite{CKS}.
  Since computations of automorphisms of quandles relates to the problem of classification of quandles, we also describe an algorithm  implemented in C for computing all quandles (up to isomorphism) of order up to nine.  \\

 In Section 2, we review the basics of quandles, give examples and describes the automorphisms  and inner automorphisms of dihedral quandles .  
 The Inner and automorphism groups of all all $seventy$ $three$ quandles of order 6 are computed in section 3.  A description of an algorithm which generates all quandles of order up to 9 (up to isomorphisms) is contained in section 4.\\
 
 \noindent
{\bf Notations} Through the paper,  the symbol $\Z_n$ will denote the set of integers modulo $n$ and $ {\Z_n}^{\times}$ will stand for the group of its units.   The dihedral group of order $2m$ will be denoted by $D_m$.   The symbol $\Sigma_n$ will stand for the symmetric group on the set $\{1,2,...,n\}$ and $A_n$ will be its alternating subgroup (even permutations).

\section{ Automorphism groups of quandles}

We start this section by reviewing the basics of quandles and give examples.

A {\it quandle}, $X$, is a set with a binary operation $(a, b) \mapsto a * b$
such that

(1) For any $a \in X$,
$a* a =a$.

(2) For any $a,b \in X$, there is a unique $x \in X$ such that 
$a= x*b$.

(3) 
For any $a,b,c \in X$, we have
$ (a*b)*c=(a*c)*(b*c). $\\
 \noindent
Axiom (2) states that for each $u \in X$, the map $S_u:X \rightarrow X$ with $S_u(x):=x*u$ is a bijection.  Its inverse will be denoted by the mapping $\overline{S}_u:X \rightarrow X$ with $\overline{S}_u(x)=x \overline{*} u$, so that $(x*u)\overline{*}u=x=(x\overline{*}u)*u.$

A {\it rack} is a set with a binary operation that satisfies 
(2) and (3).

Racks and quandles have been studied in, for example, 
\cite{FR,Joyce,Matveev}.

\noindent
The axioms for a quandle correspond respectively to the 
Reidemeister moves of type I, II, and III 
(see 
\cite{FR}, for example). 

Here are some typical examples of quandles.


\begin{list}{--}{}
\item
Any set $X$ with the operation $x*y=x$ for any $x,y \in X$ is
a quandle called the {\it trivial} quandle.
The trivial quandle of $n$ elements is denoted by $T_n$.

\item
A group $X=G$ with
$n$-fold conjugation
as the quandle operation: $a*b=b^{-n} a b^n$.

\item
Let $n$ be a positive integer.
For elements  
$i, j \in \Z_n$ (integers modulo $n$), 
define
$i\ast j \equiv 2j-i \pmod{n}$.
Then $\ast$ defines a quandle
structure  called the {\it dihedral quandle},
  $R_n$.
This set can be identified with  the
set of reflections of a regular $n$-gon
  with conjugation
as the quandle operation.
\item
Any $\Lambda (={\Z }[T, T^{-1}])$-module $M$
is a quandle with
$a*b=Ta+(1-T)b$, $a,b \in M$, called an {\it  Alexander  quandle}.
Furthermore for a positive integer
$n$, a {\it mod-$n$ Alexander  quandle}
${\Z }_n[T, T^{-1}]/(h(T))$
is a quandle
for
a Laurent polynomial $h(T)$.
The mod-$n$ Alexander quandle is finite
if the coefficients of the
highest and lowest degree terms
of $h$
  are units in $\Z_n$. 
  \end{list}


\bigskip
\noindent
A function $f: (X, *) \rightarrow  (Y,\rt)$ between quandles $X$ and $Y$  is a {\it homomorphism}
if $f(a \ast b) = f(a) \rt f(b)$ 
for any $a, b \in X$.   We will denote the group of automorphisms of the quandle $X$ by $Aut(X)$.  Axioms (2) and (3) respectively state that for each $u \in X$, the map $S_u:X \rightarrow X$ 
is respectively a bijection and a quandle homomorphism.   Lets call the subgroup of $Aut(X)$, generated by the $symmetries$ $S_x$, the $inner$ automorphism group of $X$ denoted by $Inn(X)$.   By axiom (3), the map $S:X \rightarrow Inn(X)$ sending $u$ to $S_u$ satisfies  the equation $S_z\;S_y=S_{y*z}\;S_z, \;\;\forall y,z \in X$, which can be written as $S_z\;S_y\;{S_z}^{-1}=S_{y*z}.$ Thus, if the group $Inn(X)$ is considered as a quandle with conjugation then the map $S$ becomes a quandle homomorphism.  As noted in \cite{AG} p 184, the map $S$ is not injective in general. The quandle $(X,*)$ is called $faithful$ when the map $S$ is injective.  If $(X,*)$ is $faithful$ then the center of $Inn(X)$ is trivial.\\

\subsection{Automorphism groups and Inner Automorphism groups of dihedral quandles}

Now we characterize the automorphisms of the dihedral quandles. For any non-zero element $a$ in $\Z_n$ and any $b \in \Z_n$, consider the mapping $f_{a,b}:\Z_n \rightarrow \Z_n$ sending $x $ to $ax+b $, called affine transformation over $\Z_n$.  

 \begin{theorem} \label{Aut}
Let $R_n=\Z_n$ be the dihedral quandle with the operation $i*j=2j-i \;(mod\; n)$.  Then the automorphism group $Aut(R_n)$ is isomorphic to the affine group  $\mbox{Aff}$$ (\Z_n)$.
\end{theorem}

{\it Proof.\/}
It is clear that for $a\neq 0$, the mapping $f_{a,b}$ (with  $f_{a,b}(x)=ax+b$)  is a quandle homomorphism.  It is a bijective mapping if and only if $a\in {\Z_n}^{\times}$.  Now we show that any quandle automorphism of $\Z_n$ (with the operation $x*y=2y-x$) is an affine transformation $f_{a,b}$ for some $a \in {\Z_n}^{\times}$ and  $b \in {\Z_n}$.  Let $f \in Aut(\Z_n)$, then $\forall x, y \in \Z_n, f(2y-x)=2f(y)-f(x)$.  Now consider the mapping $g:\Z_n \rightarrow \Z_n$ given by $g(x)=f(x)-f(0)$.  The mapping $g$ also satisfies $g(2y-x)=2g(y)-g(x)$.  We have $g(0)=0$ and thus $g(-a)= - g(a)$.  We now prove linearity of $g$, that is $g(\lambda x)=\lambda g(x)$ for any $\lambda \in \Z_n$.  We have $g(2b-a)=2g(b)-g(a)$, thus $g(2b)=2g(b)$ and by induction on even integers $g(2k a)= 2k g(a),$ for all $k$.  Now we do induction on odd integers:  $g[(2k+1)a]=g[2ka -(-a)]=2kg(a)-g(-a)=2kg(a)+g(a)=(2k+1)g(a).$
Now $g$ is  a bijection if and only if $\lambda \in  {\Z_n}^{\times}$ which ends the proof
\qed\\

Since the affine group $\mbox{Aff}$$ (\Z_n)$  is semi-direct product group  $\Z_n \rtimes {\Z_n}^{\times}$,  we have

\begin{corollary}
The cardinal of 
$Aut(\Z_n)$ is $n \;\phi(n)$, where $\phi$ denotes the Euler function.
\end{corollary}
For the dihedral quandle $R_n=\Z_n$ and for each $i \in \Z_n$ the  symmetry $S_i$ given by $S_i(j)= 2i-j \;(mod\; n)$, can be though of as a reflection of a regular $n$-gon.  If n is odd, the axis of symmetry of $S_i$ connects the vertex $i$ to the mid-point of the side opposite to $i$. If $n=2m$ is even, the axis of symmetry of $S_i$ passes through the opposite vertices $i$ and $i+m$ (mod $2m$).  From these observations, we have the easy characterization of the inner automorphism group of dihedral quandles given by the following

 \begin{theorem} \label{Aut}
The inner automorphism group $Inn(R_n)$ of the dihedral quandle $R_n$  is isomorphic to the dihedral group  $D_{\frac{m}{2}}$ of order $m$ where $m$ is the least common multiple of $n$ and $2$.
\end{theorem}

 \begin{theorem} \label{InnerConj}
 Let $G$ be a group and let the quandle $X$ be the group $G$ as a set with the conjugation $x*y=yxy^{-1}$ as operation.  This quandle is usually denoted by $Conj(G)$. Then the Inner automorphism group of $X$ is isomorphic (as a group) to the quotient of $G$ by its center $Z(G)$.  
  \end{theorem}
  
{\it Proof.\/}
  The proof is straightforward from the fact that in this case the surjective map $S: X \rightarrow Inn(X)$ sending $a \in X$ to $S_a$ is a quandle homomorphism with kernel the center $Z(G)$ of $G$. \qed\\
 
 \noindent 
  {\bf Example} The symetric group $\Sigma_3$ is the smallest group with trivial center then $Inn(Conj(\Sigma_3))\cong \Sigma_3.$\\
  
\noindent
The converse of theorem  \ref{InnerConj} is also true, namely if $(X,*)$ is a quandle for which the map $S: X \rightarrow Inn(X)$ is one-to-one and onto then $(X,*) \cong Conj(Inn(X))$ with $Z(Inn(X))$ being trivial group.\\
\noindent
An interesting question would be to calculate the automorphism groups $Aut(Conj(G))$.  Obviously for the symmetric group $\Sigma_3$, we have $Aut(Conj(\Sigma_3))\cong Inn(Conj(\Sigma_3))\cong \Sigma_3. $\\

\section{Automorphism and Inner Automorphism groups of quandles of order $6$} 
In this section, we compute the automorphism groups and the inner automorphism groups of all $seventy$ $three$ quandles of order six.   The computation is accomplished with the help of the software Maple which also allows the computation of the  inner and automorphism groups for quandles of order $7$ and $8$.  Since the numbers of isomorphism classes of quandles of order $7$ and $8$ are respectively $298$ and $1581$, we decided not to include these two cases in this paper. \\
We describe each quandle $Q_j$ of order $6$ for $1 \leq j \leq 73$ by explicitly giving each symmetry $S_k$ for $1 \leq k \leq 6$, in terms of products of disjoint cycles.  The symmetries are the columns in the Cayley table of the quandle.  For example  the quandle, denoted $Q_{46}$ in table $2$  below, with the Cayley table

$$
\left[ \begin{array}{rrrrrr}
1 &1&1&1&1&1\\
2 &2&5&5&2&5\\ 
3 &4&3&3&4&4\\ 
4 &3&4&4&3&3\\ 
5 &5&2&2&5&2\\ 
6 &6&6&6&6&6
\end{array} \right]
$$
is described by the permutations of the six elements set $\{1,2,3,4,5,6\}$,   $S_1 = (1)$,\; $S_2 = (34)$,\; $S_3 =(25)$,\; $S_4 = (25)$,\; $S_5 =(34)$,\;$S_6 = (25)(34)$.  Here and through the rest of the paper, every permutation is written as a product of transpositions.  For example, $S_1=(1)$ means that $S_1$ is the identity permutation.   The permutation $S_4 = (25)$ stands for the transposition sending $2$ to $5$ and $S_6 = (25)(34)$ stands for the product of the two transpositions $(25)$ and $(34)$.  \\
In this example $Aut(Q_{46})=D_4,$ the dihedral group of $8$ elements and $Inn(Q_{46})=\Z_2 \times \Z_2$ is the direct product of two copies of $\Z_2$. Another example given in table \ref{tableAut} is $Aut(Q_{49})=D_5$ the dihedral group of order $10$ and $Inn(Q_{46})=\Z_5 \rtimes \Z_4,$ the semidirect product of the cyclic group $\Z_5$ by  $\Z_4$.  
\begin{table}
\begin{center}
\begin{tabular}{||l|c|c||} \hline \hline
Quandle & Disjoint Cycle Notation for the Columns of the Quandle  
\\ \hline

$Q_1$ & $(1),(1),(1),(1),(1),(1)$   \\ \hline
$Q_2$ & $(1),(1),(1),(1),(1),(12)$   \\ \hline
$Q_3$ & $(1),(1),(1),(1),(1),(132)$   \\ \hline
$Q_4$ & $(1),(1),(1),(1),(1),(1243)$   \\ \hline
$Q_5$ & $(1),(1),(1),(1),(1),(12)(34)$     \\ \hline
$Q_6$ & $(1),(1),(1),(1),(1),(15234)$   \\ \hline
$Q_7$ & $(1),(1),(1),(1),(1),(134) (25)$   \\ \hline
$Q_8$ & $(1),(1),(1),(1),(12),(12)$       \\ \hline
$Q_9$ & $(1),(1),(1),(1),(12),(12) (34)$      \\ \hline
$Q_{10}$ & $(1),(1),(1),(1),(12),(34)$   \\ \hline
$Q_{11}$ & $(1),(1),(1),(1),(132),(132)$    \\ \hline
$Q_{12}$ & $(1),(1),(1),(1),(132),(123)$     \\ \hline
$Q_{13}$ & $(1),(1),(1),(1),(1243),(1243)$      \\ \hline
$Q_{14}$ & $(1),(1),(1),(1),(1243),(1342)$     \\ \hline
$Q_{15}$ & $(1),(1),(1),(1),(1243),(14)(23)$        \\ \hline
$Q_{16}$ & $(1),(1),(1),(1),(12)(34),(12)(34)$        \\ \hline
$Q_{17}$ & $(1),(1),(1),(1),(12)(34),(13)(24)$     \\ \hline
$Q_{18}$ & $(1),(1),(1),(12),(12),(12)$   \\ \hline
$Q_{19}$ & $(1),(1),(1),(12),(12),(12)(45)$      \\ \hline
$Q_{20}$ & $(1),(1),(1),(12),(12),(45)$    \\ \hline
$Q_{21}$ & $(1),(1),(1),(132),(132),(132)$    \\ \hline
$Q_{22}$ & $(1),(1),(1),(132),(132),(123)$      \\ \hline
$Q_{23}$ & $(1),(1),(1),(132),(132),(45)$      \\ \hline
$Q_{24}$ & $(1),(1),(1),(132),(132),(123)(45)$    \\ \hline
$Q_{25}$ & $(1),(1),(1),(132),(132),(132)(45)$  \\ \hline
$Q_{26}$ & $(1),(1),(1),(12)(56),(12)(46),(12)(45)$  \\ \hline
$Q_{27}$ & $(1),(1),(1),(12)(56),(13)(46),(23)(45)$        \\ \hline
$Q_{28}$ & $(1),(1),(1),(56),(46),(45)$     \\ \hline
$Q_{29}$ & $(1),(1),(1),(123)(56),(123)(46),(123)(45)$   \\ \hline
$Q_{30}$ & $(1),(1),(12),(12),(12),(12)$     \\ \hline
$Q_{31}$ & $(1),(1),(12),(12),(12),(12)(34)$   \\ \hline
$Q_{32}$ & $(1),(1),(12),(12),(12),(34)$     \\ \hline
$Q_{33}$ & $(1),(1),(12),(12),(12),(345)$    \\ \hline
$Q_{34}$ & $(1),(1),(12),(12),(12),(12)(345)$      \\ \hline
$Q_{35}$ & $(1),(1),(12),(12),(12)(34),(12)(34)$       \\ \hline
$Q_{36}$ & $(1),(1),(12),(12),(12)(34),(34)$    \\ \hline
$Q_{37}$ & $(1),(1),(12),(12),(34),(34)$   \\ \hline

\end{tabular}\end{center}
\caption{Quandles of order 6 in term of disjoint cycles of columns - part 1}
\label{qmoduletable3}
\end{table}


\begin{table}
\begin{center}
\begin{tabular}{||l|c|c||} \hline \hline
Quandle  & Disjoint Cycle Notation for the Columns of the Quandle  
\\ \hline

$Q_{38}$ & $(1),(1),(12),(12)(56),(12)(46),(12)(45) $     \\ \hline
$Q_{39}$ & $(1),(1),(12),(56),(46),(45)$     \\ \hline
$Q_{40}$ & $(1),(1),(12)(45),(12)(36),(12)(36),(12)(45)$   \\ \hline
$Q_{41}$ & $(1),(1),(12)(45),(36),(36),(12)(45)$   \\ \hline
$Q_{42}$ & $(1),(1),(45),(36),(36),(45)$     \\ \hline
$Q_{43}$ & $(1),(1),(456),(365),(346),(354)$   \\ \hline
$Q_{44}$ & $(1),(1),(12)(456),(12)(365),(12)(346),(12)(354)$    \\ \hline
$Q_{45}$ & $(1),(34), (25), (25),(34),(34)$       \\ \hline
$Q_{46}$ & $(1),(34),(25),(25),(34),(25)(34)$     \\ \hline
$Q_{47}$ & $(1),(34),(256),(256),(34),(34)$   \\ \hline
$Q_{48}$ & $(1),(354),(26)(45),(26)(35),(26)(34),(345)$   \\ \hline
$Q_{49}$ & $(1),(36)(45),(25)(46),(23)(56),(26)(34),(24)(35)$     \\ \hline
$Q_{50}$ & $(1),(3546),(2456),(2365),(2643),(2534)$      \\ \hline
$Q_{51}$ & $(1),(3546),(2564),(2653),(2436),(2345)$    \\ \hline
$Q_{52}$ & $(23),(13),(12),(56),(46),(45)$      \\ \hline
$Q_{53}$ & $(23),(14),(14),(23),(23),(23)$       \\ \hline
$Q_{54}$ & $(23),(14),(14),(23),(23),(14)(23)$      \\ \hline
$Q_{55}$ & $(23),(14),(14),(23),(23),(14)$    \\ \hline
$Q_{56}$ & $(23),(14),(14),(23),(14)(23),(14)(23)$      \\ \hline
$Q_{57}$ & $(23),(154),(154),(23),(23),(23)$    \\ \hline
$Q_{58}$ & $(23),(154),(154),(23),(23),(154)(23)$    \\ \hline
$Q_{59}$ & $(23),(154),(154),(23),(23),(154)$       \\ \hline
$Q_{60}$ & $(23),(154),(154),(23),(23),(145)$    \\ \hline
$Q_{61}$ & $(23),(154),(154),(23),(23),(145)(23)$      \\ \hline
$Q_{62}$ & $(23),(45),(45),(16)(23),(16)(23),(23)$   \\ \hline
$Q_{63}$ & $(23),(45),(45),(16),(16),(23)$    \\ \hline
$Q_{64}$ & $(23),(1564),(1564),(23),(23),(23)$        \\ \hline
$Q_{65}$ & $(23),(15)(46),(15)(46),(23),(23),(23)$    \\ \hline
$Q_{66}$ & $(23),(15)(46),(15)(46),(15)(23),(23),(15)(23)$  \\ \hline
$Q_{67}$ & $(243),(165),(165),(165),(243),(243)$          \\ \hline
$Q_{68}$ & $(2354),(1463),(1265),(1562),(1364),(2453)$  \\ \hline
$Q_{69}$ & $(2354),(16)(34),(16)(25),(16)(25),(16)(34),(2453)$   \\ \hline
$Q_{70}$ & $(23)(45),(15)(36),(14)(26),(15)(36),(14)(26),(23)(45)$      \\ \hline
$Q_{71}$ & $(23)(45),(15)(46),(14)(56),(16)(23),(16)(23),(23)(45)$     \\ \hline
$Q_{72}$ & $(23)(45),(13)(46),(12)(56),(15)(26),(14)(36),(24)(35)$         \\ \hline
$Q_{73}$ & $(23)(45),(16)(45),(16)(45),(16)(23),(16)(23),(23)(45)$    \\ \hline

\end{tabular}\end{center}
\caption{Quandles of order 6 in term of disjoint cycles of columns - part 2}
\label{qmoduletable3}
\end{table}


\begin{table}\label{tableAut}
\begin{center}
\begin{tabular}{||l||c|c||c||c|c|c||} \hline \hline
Quandle X & Inn(X) & Aut(X)    & Quandle X & Inn(X) & Aut(X)   
\\ \hline 

$Q_1$ & \{1\}  &  $\Sigma_6$ & $Q_{38}$ & $D_3 \times \Z_2$ & $D_3 \times \Z_2$       \\ \hline
$Q_2$ & $\Z_2 $ & $D_3 \times \Z_2$   &$Q_{39}$  &$D_3 \times \Z_2$  &$D_3 \times \Z_2$      \\ \hline
$Q_3$ & $\Z_3 $ & $\Z_6$   &$Q_{40}$  & $\Z_2 \times \Z_2$ & $D_4 \times \Z_2$     \\ \hline
$Q_4$ & $\Z_4 $ & $\Z_4$   &$Q_{41}$  & $\Z_2 \times \Z_2 $  &  $\Z_2 \times \Z_2 \times \Z_2$     \\ \hline
$Q_5$ & $\Z_2 $ & $D_4$   &$Q_{42}$  & $\Z_2 \times \Z_2$  &  $D_4 \times \Z_2$      \\ \hline
$Q_6$ & $\Z_5 $ & $\Z_5$   &$Q_{43}$  & $A_4$  & $A_4\times \Z_2$      \\ \hline
$Q_7$ & $\Z_6 $ & $\Z_6$   &$Q_{44}$  & $A_4\times \Z_2$  & $A_4\times \Z_2$      \\ \hline
$Q_8$ & $\Z_2 $ & $\Z_2 \times \Z_2 \times \Z_2$   &$Q_{45}$  &$\Z_2 \times \Z_2 $  & $\Z_2 \times \Z_2 $     \\ \hline
$Q_9$ & $\Z_2 \times \Z_2 $ & $\Z_2 \times \Z_2$   &$Q_{46}$  &$\Z_2 \times \Z_2 $  & $D_4$     \\ \hline
$Q_{10}$ & $\Z_2 \times \Z_2$ & $D_4 $   &$Q_{47}$  &$\Z_6 $  & $\Z_6$     \\ \hline
$Q_{11}$ & $\Z_3  $ & $\Z_6$   &$Q_{48}$  &$D_3$  &$D_3$      \\ \hline
$Q_{12}$ & $\Z_3  $ & $D_3 $   &$Q_{49}$  &$D_5$  & $\Z_5\rtimes\Z_4$    \\ \hline
$Q_{13}$ & $\Z_4  $ & $\Z_4 \times \Z_2$   &$Q_{50}$  & $\Z_5\rtimes\Z_4$  & $\Z_5\rtimes\Z_4$      \\ \hline
$Q_{14}$ & $\Z_4 $ & $D_4$   &$Q_{51}$  & $\Z_5\rtimes\Z_4$  &  $\Z_5\rtimes\Z_4$     \\ \hline
$Q_{15}$ & $\Z_4  $ & $\Z_4 $   &$Q_{52}$  &$D_3 \times\ D_3$   & $(D_3 \times\ D_3) \rtimes\Z_2$     \\ \hline
$Q_{16}$ & $\Z_2  $ & $D_4 \times \Z_2$   &$Q_{53}$  &$\Z_2 \times \Z_2 $  &$\Z_2 \times \Z_2 \times \Z_2$      \\ \hline
$Q_{17}$ & $\Z_2 \times \Z_2 $ & $D_4$   &$Q_{54}$  & $\Z_2 \times \Z_2 $ & $\Z_2 \times \Z_2 $     \\ \hline
$Q_{18}$ & $\Z_2 $ & $D_3 \times \Z_2$   &$Q_{55}$  &$\Z_2 \times \Z_2 $  & $D_4$      \\ \hline
$Q_{19}$ & $\Z_2 \times \Z_2 $ & $\Z_2 \times \Z_2$   &$Q_{56}$  & $\Z_2 \times \Z_2 $ &  $D_4 \times \Z_2$    \\ \hline
$Q_{20}$ & $\Z_2 \times \Z_2 $ & $\Z_2 \times \Z_2$   &$Q_{57}$  &$\Z_6 $  &$\Z_6 $      \\ \hline
$Q_{21}$ & $\Z_3 $ & $D_3 \times \Z_3$   &$Q_{58}$  &$\Z_6 $  & $\Z_6 $     \\ \hline
$Q_{22}$ & $\Z_3 $ & $\Z_6$   &$Q_{59}$  &$\Z_6 $  & $\Z_6 $     \\ \hline
$Q_{23}$ & $\Z_6 $ & $\Z_6$   &$Q_{60}$  &$\Z_6 $  &$\Z_6 $      \\ \hline
$Q_{24}$ & $\Z_6 $ & $\Z_6$   &$Q_{61}$  &$\Z_6 $  &$\Z_6 $      \\ \hline
$Q_{25}$ & $\Z_6 $ & $\Z_6$   &$Q_{62}$  & $\Z_2 \times \Z_2 \times \Z_2$  &      $\Z_2 \times \Z_2 \times \Z_2$ \\ \hline
$Q_{26}$ & $D_3 $ & $D_3 \times \Z_2$   &$Q_{63}$  & $\Z_2 \times \Z_2 \times \Z_2$  &  $A_4\times \Z_2$    \\ \hline
$Q_{27}$ & $D_3 $ & $D_3 $   &$Q_{64}$  &$\Z_4 \times \Z_2 $  &$\Z_4 \times \Z_2 $      \\ \hline
$Q_{28}$ & $D_3 $ & $D_3 \times\ D_3$   &$Q_{65}$  &$\Z_2 \times \Z_2 $  & $D_4 \times \Z_2$  \\ \hline
$Q_{29}$ & $D_3 \times \Z_3 $ & $D_3 \times \Z_3$   &$Q_{66}$  &  $\Z_2 \times \Z_2 \times \Z_2$ &      $\Z_2 \times \Z_2 \times \Z_2$ \\ \hline
$Q_{30}$ & $\Z_2 $ & $\Sigma_4 \times \Z_2$   &$Q_{67}$  & $\Z_3 \times \Z_3 $  & $D_3 \times \Z_3$      \\ \hline
$Q_{31}$ & $\Z_2 \times \Z_2 $ & $\Z_2 \times \Z_2$   &$Q_{68}$  & $\Sigma_4$  &  $\Sigma_4$     \\ \hline
$Q_{32}$ & $\Z_2 \times \Z_2 $ & $\Z_2 \times \Z_2$   &$Q_{69}$  & $D_4 $  &  $D_4 $     \\ \hline
$Q_{33}$ & $\Z_6 $ & $\Z_6$   &$Q_{70}$  & $D_3 $  & $D_3 \times \Z_2$     \\ \hline
$Q_{34}$ & $\Z_6 $ & $\Z_6$   &$Q_{71}$  & $D_4 $  & $D_4 $      \\ \hline
$Q_{35}$ & $\Z_2 \times \Z_2 $ & $\Z_2 \times \Z_2 \times \Z_2$   &$Q_{72}$  & $\Sigma_4$ &$\Sigma_4$      \\ \hline
$Q_{36}$ & $\Z_2 \times \Z_2 $ & $\Z_2 \times \Z_2$   &$Q_{73}$  &$\Z_2 \times \Z_2 $  & $\Sigma_4 \times \Z_2 $     \\ \hline
$Q_{37}$ & $\Z_2 \times \Z_2 $ & $\Z_2 \times \Z_2 \times \Z_2$   &  &  &      \\ \hline
\hline
\end{tabular}\end{center}
\caption{A table of the quandles  of order 6 with their Inner and Automorphism groups}
\label{qmoduletable3}
\end{table}

\newpage

\section{Algorithm description}

In the quest of finding computationally the quandles of certain order up to
isomorphism, we are cursed by the fact that any sort of naive algorithm will
take an exponential time (in the order of the quandle) to do such task.
Therefore, we are required to exploit structural or logical aspects of the
quandle theory to reduce the running time at least by a proportional factor,
making the algorithm `less-galactic', in CS jargon.

\subsection{Phase 1:\ List generation.}

In initial versions of the \emph{quandles algorithm} \cite{HMN}, the set of
all matrices such that every row is a permutation $\left[  n\right]  $, is
generated. After this, the matrices that do not correspond to the operation
table of a quandle (i.e., such that do not satisfy the quandle axiom), are
ruled out. We call this initial process the \emph{list generation}, and its
purpose is to \textbf{list a set of quandles such that among then, we are
guaranteed to find representatives for all isomorphic classes of quandles of
order }$n$. A further improvement in this process consists in verifying the
quandles axiom \emph{on-line}, this means that, during the generation of the
matrices, the axioms are immediately verified, a process that was also carried
out in \cite{HMN}. We elaborate this improvement to a higer level:\ Besides
verifying the quandle axioms on-line, we also fill in online, entries that are
implied by the quandle axioms. To exemplify such process, suppose that at a
certain step our algorithm has completed the following partial table of a
quandle
\[%
\begin{bmatrix}
a & a &  &  & \\
c & b & b &  & \\
b & c & c &  & \\
&  &  & d & \\
&  &  &  & e
\end{bmatrix}
\]
then, by use of the property $\left(  a\ast c\right)  \ast a=a\ast\left(
c\ast a\right)  $, we have that $\left(  a\ast c\right)  \ast a=a$. Therefore
$\left(  a\ast c\right)  =a\bar{\ast}a=a$, so that the table completes as%
\[%
\begin{bmatrix}
a & a & a &  & \\
c & b & b &  & \\
b & c & c &  & \\
&  &  & d & \\
&  &  &  & e
\end{bmatrix}
.
\]
A more interesting example is the following. Starting with the following
partial quandle table
\[%
\begin{bmatrix}
a & a & a & b & \\
c & b & b &  & \\
b & c & c &  & \\
&  &  & d & \\
&  &  & c & e
\end{bmatrix}
,
\]

through the application of the quandle axioms several times, we complete some
fewer entries, concluding at the end that such partial table cannot be
extended to a valid quandle table:
\begin{align*}
& \overset{\left(  e\ast a\right)  \ast d=\left(  e\ast d\right)  \ast\left(
a\ast d\right)  }{\rightarrow}%
\begin{bmatrix}
a & a & a & b & \\
c & b & b &  & \\
b & c & c &  & \\
&  &  & d & \\
e &  &  & c & e
\end{bmatrix}
\overset{\left(  a\ast e\right)  \ast d=\left(  a\ast d\right)  \ast\left(
e\ast d\right)  }{\rightarrow}%
\begin{bmatrix}
a & a & a & b & a\\
c & b & b &  & \\
b & c & c &  & \\
&  &  & d & \\
e &  &  & c & e
\end{bmatrix}
\\
& \overset{\text{uniqueness}}{\rightarrow}%
\begin{bmatrix}
a & a & a & b & a\\
c & b & b &  & \\
b & c & c &  & \\
d &  &  & d & \\
e &  &  & c & e
\end{bmatrix}
\end{align*}
and this last table contradicts the axiom $\left(  a\ast d\right)  \ast
a=a\ast\left(  d\ast a\right)  $.

In general, the rules that are used for this `completion' process are the following:

Suppose that $j\ast i=k$, then

\emph{Rule 1: }$k\ast a=\left(  j\ast a\right)  \ast\left(  i\ast a\right)  $

\begin{enumerate}
\item If $\left(  j\ast a\right)  \ast\left(  i\ast a\right)  $ cannot be
retrieved from the table and $k\ast a$, $i\ast a$ and $j\ast a$ can be
retrieved from the table, then

\begin{enumerate}
\item If $\left(  k\ast a\right)  \bar{\ast}\left(  i\ast a\right)  $ can be
retrieved from the table, the table is not valid.

\item Otherwise, necesarily $\left(  j\ast a\right)  \ast\left(  i\ast
a\right)  =k\ast a$.\newline
\end{enumerate}

\item Otherwise, if $k\ast a$ cannot be retrieved from the table and $\left(
j\ast a\right)  \ast\left(  i\ast a\right)  $ can be retrieved from the table, then

\begin{enumerate}
\item If $\left(  \left(  j\ast a\right)  \ast\left(  i\ast a\right)  \right)
\bar{\ast}a$ can be retrieved from the table, the table is not valid.

\item Otherwise, necesarily $k\ast a=\left(  j\ast a\right)  \ast\left(  i\ast
a\right)  $.
\end{enumerate}

\item Otherwise, if $\left(  j\ast a\right)  \ast\left(  i\ast a\right)  $ and
$k\ast a$ can be retrieved from the table and $\left(  j\ast a\right)
\ast\left(  i\ast a\right)  \neq k\ast a$, the table is not valid.
\end{enumerate}

\emph{Rule 2: }$\left(  a\ast j\right)  \ast i=\left(  a\ast i\right)  \ast k$

\begin{enumerate}
\item If $\left(  a\ast i\right)  \ast k$ cannot be retrieved from the table
and $\left(  a\ast j\right)  \ast i$ and $a\ast i$ can be retrieved from the
table, then

\begin{enumerate}
\item If $\left(  \left(  a\ast j\right)  \ast i\right)  \bar{\ast}k$ can be
retrieved from the table, the table is not valid.

\item Otherwise, necesarily $\left(  a\ast i\right)  \ast k=$\emph{\ }$\left(
a\ast j\right)  \ast i$.\newline
\end{enumerate}

\item Otherwise, if $\left(  a\ast j\right)  \ast i$ cannot be retrieved from
the table and $\left(  a\ast j\right)  $ and $\left(  a\ast i\right)  \ast k$
can be retrieved from the table, then

\begin{enumerate}
\item If $\left(  \left(  a\ast i\right)  \ast k\right)  \bar{\ast}i$ can be
retrieved from the table, the table is not valid.

\item Otherwise, necesarily $\left(  a\ast j\right)  \ast i=\left(  a\ast
i\right)  \ast k$.
\end{enumerate}

\item Otherwise, if $\left(  a\ast j\right)  \ast i$ and $\left(  a\ast
i\right)  \ast k$ can be retrieved from the table and $\left(  a\ast j\right)
\ast i\neq\left(  a\ast i\right)  \ast k$, the table is not valid.
\end{enumerate}

\emph{Rule 3: }$\left(  j\ast a\right)  \ast i=k\ast\left(  a\ast i\right)  $

\begin{enumerate}
\item If $k\ast\left(  a\ast i\right)  $ cannot be retrieved from the table
and $\left(  j\ast a\right)  \ast i$ and $a\ast i$ can be retrieved from the
table, then

\begin{enumerate}
\item If $\left(  \left(  j\ast a\right)  \ast i\right)  \bar{\ast}\left(
a\ast i\right)  $ can be retrieved from the table, the table is not valid.

\item Otherwise, necesarily $k\ast\left(  a\ast i\right)  =$\emph{\ }$\left(
j\ast a\right)  \ast i$.\newline
\end{enumerate}

\item Otherwise, if $\left(  j\ast a\right)  \ast i$ cannot be retrieved from
the table and $\left(  j\ast a\right)  $ and $k\ast\left(  a\ast i\right)  $
can be retrieved from the table, then

\begin{enumerate}
\item If $\left(  k\ast\left(  a\ast i\right)  \right)  \bar{\ast}i$ can be
retrieved from the table, the table is not valid.

\item Otherwise, necesarily $\left(  j\ast a\right)  \ast i=k\ast\left(  a\ast
i\right)  $.
\end{enumerate}

\item Otherwise, if \emph{\ }$\left(  j\ast a\right)  \ast i$ and
$k\ast\left(  a\ast i\right)  $ can be retrieved from the table and
\emph{\ }$\left(  j\ast a\right)  \ast i\neq k\ast\left(  a\ast i\right)  $,
the table is not valid.
\end{enumerate}

\emph{Rule 4: }$\left(  \left(  j\bar{\ast}a\right)  \ast\left(  i\bar{\ast
}a\right)  \right)  \ast a=k$

\begin{enumerate}
\item If $\left(  \left(  j\bar{\ast}a\right)  \ast\left(  i\bar{\ast
}a\right)  \right)  \ast a$ cannot be retrieved from the table and $\left(
\left(  j\bar{\ast}a\right)  \ast\left(  i\bar{\ast}a\right)  \right)  $ can
be retrieved from the table, then

\begin{enumerate}
\item If $k\bar{\ast}a$ can be retrieved from the table, the table is not valid.

\item Otherwise, necesarily $\left(  \left(  j\bar{\ast}a\right)  \ast\left(
i\bar{\ast}a\right)  \right)  \ast a=k$.\newline
\end{enumerate}

\item Otherwise, if $\left(  \left(  j\bar{\ast}a\right)  \ast\left(
i\bar{\ast}a\right)  \right)  \ast a$ can be retrieved from the table and
$\left(  \left(  j\bar{\ast}a\right)  \ast\left(  i\bar{\ast}a\right)
\right)  \ast a\neq k$, the table is not valid.
\end{enumerate}

\emph{Rule 5: }$k=\left(  \left(  j\bar{\ast}a\right)  \ast i\right)
\ast\left(  a\ast i\right)  $

\begin{enumerate}
\item If $\left(  \left(  j\bar{\ast}a\right)  \ast i\right)  \ast\left(
a\ast i\right)  $ cannot be retrieved from the table and $\left(  a\ast
i\right)  $ and\ $\left(  \left(  j\bar{\ast}a\right)  \ast i\right)  $ can be
retrieved from the table, then

\begin{enumerate}
\item If $k\bar{\ast}\left(  a\ast i\right)  $ can be retrieved from the
table, the table is not valid.

\item Otherwise, necesarily $\left(  \left(  j\bar{\ast}a\right)  \ast
i\right)  \ast\left(  a\ast i\right)  =k$.\newline
\end{enumerate}

\item Otherwise, if $\left(  \left(  j\bar{\ast}a\right)  \ast i\right)
\ast\left(  a\ast i\right)  $ can be retrieved from the table and
$k\neq\left(  \left(  j\bar{\ast}a\right)  \ast i\right)  \ast\left(  a\ast
i\right)  $, the table is not valid.
\end{enumerate}

Another easy improvement, which certainly reduces considerably the size of the
list of quandles to output in this first step of the quandles algorithm, comes
from elementary logic: When you are trying to generate all the models of
cardinality $n$ of a theory (in our case the theory of quandles), we can start
introducing constants and the corresponding relations between these constants
one by one (in a valid way), until we get $n$ constants (so, the possible ways
to generate the relations between constants will correspond to the models of
the theory). This is exactly what any algorithm will do, just in the language
of logic, but the point to emphasize is that, when a new constant is
introduced, the name of such constant is irrelevant. This is a trivial logic
fact, but one that was not used in previous versions of this listing
procedure. For example, if we aim to complete the entry $b\ast a$ of the
partial table
\[%
\begin{bmatrix}
a &  &  &  & \\
& b &  &  & \\
&  & c &  & \\
&  &  & d & \\
&  &  &  & e
\end{bmatrix}
\text{,}%
\]
then among the options $b\ast a=c$, $b\ast a=d$ and $b\ast a=e$, the choice is
irrelevant, because at such step, the constants $c,d,e$ are not in context.

The following are some benchmarks concerning this first step of the process:%

\[%
\begin{bmatrix}
\text{size} & \text{quandles} & \text{time (sec.)}\\
2 & 1 & 0\\
3 & 5 & 0\\
4 & 27 & 0\\
5 & 190 & 0\\
6 & 1833 & 0\\
7 & 22104 & 1\text{ to }2\\
8 & 359859 & 24\text{ to }34
\end{bmatrix}
\]

\subsection{Phase 2:\ Isomorphic comparison}

After the previous listing procedure has been elaborated (or more precisely,
\textbf{while} the listing procedure is elaborated), we want to eliminate
irrelevant quandles, that is, we want to leave only one representative per
isomorphism class. For such comparison process, instead of doing a brute force
algorithm that takes all possible bijections and checks for isomorphic
equivalence, we can do two things:

(1) Use simple invariant checks, like number of cycles in every row action, to
discard rapidly some nonisomorphic pairs of quandles.

(2) Use the quandle axioms to reduce the complexity of the isomorphic
comparison process.

Regarding (2), we employ the quandle axioms to extend appropriately a partial
isomorphism among valid possibilities, using the following rules:

Suppose that $\phi\left(  i\right)  =j$.

\emph{Rule 1: }$\phi\left(  b\right)  \ast^{\prime}j=\phi\left(  b\ast
i\right)  $

\begin{enumerate}
\item If $\phi\left(  b\right)  $ and $\phi\left(  b\ast i\right)  $ are
defined, and $\phi\left(  b\right)  \ast^{\prime}j\neq\phi\left(  b\ast
i\right)  $, then the isomorphism is not valid.

\item If $\phi\left(  b\ast i\right)  $ is not defined and $\phi\left(
b\right)  $ is defined, necessarily $\phi\left(  b\ast i\right)  =\phi\left(
b\right)  \ast^{\prime}j$, and this may or may not contradict the injectivity
of $\phi$.

\item If $\phi\left(  b\ast i\right)  $ is defined and $\phi\left(  b\right)
$ is not defined, necessarily $\phi\left(  b\right)  =\phi\left(  b\ast
i\right)  \bar{\ast}^{\prime}j$, and this may or may not contradict the
injectivity of $\phi$.
\end{enumerate}

\emph{Rule 2: }$j\ast^{\prime}\phi\left(  b\right)  =\phi\left(  i\ast
b\right)  $

\begin{enumerate}
\item If $\phi\left(  b\right)  $ and $\phi\left(  i\ast b\right)  $ are
defined, and $j\ast^{\prime}\phi\left(  b\right)  \neq\phi\left(  i\ast
b\right)  $, then the isomorphism is not valid.

\item If $\phi\left(  i\ast b\right)  $ is not defined and $\phi\left(
b\right)  $ is defined, necessarily $\phi\left(  i\ast b\right)  =j\ast
^{\prime}\phi\left(  b\right)  $, and this may or may not contradict the
injectivity of $\phi$.
\end{enumerate}

\emph{Rule 3: }$\phi\left(  i\bar{\ast}b\right)  \ast^{\prime}\phi\left(
b\right)  =j$

\begin{enumerate}
\item If $\phi\left(  b\right)  $ and $\phi\left(  i\bar{\ast}b\right)  $ are
defined, and $\phi\left(  i\bar{\ast}b\right)  \ast^{\prime}\phi\left(
b\right)  \neq j$, then the isomorphism is not valid.

\item If $\phi\left(  b\right)  $ is defined and $\phi\left(  i\bar{\ast
}b\right)  $ is not defined, necessarily, $\phi\left(  i\bar{\ast}b\right)
=j\bar{\ast}^{\prime}\phi\left(  b\right)  $, and this may or may not
contradict the injectivity of $\phi$.
\end{enumerate}

For the following benchmark, we do an exhaustive algorithm for isomorphism
comparison. Notice that is tractable up to $n=6$.%

\[%
\begin{bmatrix}
\text{size} & \text{quandles} & \text{time (sec.)}\\
2 & 1 & 0\\
3 & 3 & 0\\
4 & 7 & 0\\
5 & 22 & 0\\
6 & 73 & 29-32
\end{bmatrix}
\]

For the following benchmark we apply the improved isomorphism comparison, by
using the rules described previously. This improves the running time by a
factor of $10$ approx.%

\[%
\begin{bmatrix}
\text{size} & \text{quandles} & \text{time (sec.)}\\
2 & 1 & 0\\
3 & 3 & 0\\
4 & 7 & 0\\
5 & 22 & 0\\
6 & 73 & 3\\
7 & 298 & 330
\end{bmatrix}
\]

\bigskip\bigskip

\emph{Checking invariants:}

\medskip

Certainly, it is not necessary to do an isomorphism comparison (improved or
not), if we know before hand that the quandles to be compared are `too
different'. Therefore, a pre-comparison of some invariants fast to calculate,
would boost the running time. For early versions of the algorithm, we
introduced invariants based on the permutation structure of the columns of the
quandle table. For example, for the following benchmark, we simply count the
total number of cycles among all columns of the quandle table. The comparison
of such invariant improves the running time by another factor of $10$:%

\[%
\begin{bmatrix}
\text{size} & \text{quandles} & \text{time (sec.)}\\
2 & 1 & 0\\
3 & 3 & 0\\
4 & 7 & 0\\
5 & 22 & 0\\
6 & 73 & 0\\
7 & 298 & 34
\end{bmatrix}
\]

For the following benchmark we go down one more level, now taking as invariant
the superset consisting of the number of cycles of every columns. This improves the running time by a factor of $4$ approx.%

\[%
\begin{bmatrix}
\text{size} & \text{quandles} & \text{time (sec.)}\\
2 & 1 & 0\\
3 & 3 & 0\\
4 & 7 & 0\\
5 & 22 & 0\\
6 & 73 & 0\\
7 & 298 & 9
\end{bmatrix}
\]

For the following benchmark we refine the previous invariant, by considering
the superset of supersets of cycle lengths of every column, At this level of
improvement, the case $n=8$ is computationally tractable.%

\[%
\begin{bmatrix}
\text{size} & \text{quandles} & \text{time (sec.)}\\
2 & 1 & 0\\
3 & 3 & 0\\
4 & 7 & 0\\
5 & 22 & 0\\
6 & 73 & 0\\
7 & 298 & 6\\
8 & 1581 & 458
\end{bmatrix}
\]
\noindent

Further improvements will be
introduced in next versions of the algorithm, whose source is available at the
web address http://people.math.gatech.edu/\symbol{126}restrepo/quandles.html.\\
Another invariants suggested by Professor Edwin Clark, which according to his experiments seem to distinguish isomorphic classes effectively, take in account the structure of the rows of the quandle.\\
The number of isomorphism of quandles of order 3, 4, 5, 6, 7, 8 and 9 we obtain are respectively 3, 7, 22, 73, 298, 1581, 11079.  These same numbers are obtained by  James McCarron in \cite{OEIS}.

\noindent
{\large\bf Acknowledgments} \
The authors would like to thank professor Edwin Clark for his help and fruitful suggestions.

\end{document}